\documentclass[12pt]{article}
\usepackage[utf8]{inputenc}
\usepackage[T1]{fontenc}
\usepackage{amsmath,amssymb,amsfonts,amsthm}
\usepackage{natbib}
\usepackage{float}
\usepackage{booktabs,multirow}
\usepackage{tabularx}
\usepackage{array}
\usepackage{geometry}
\usepackage[toc,page]{appendix}
\usepackage{xcolor}
\newtheorem{theorem}{Theorem}%[section]
\newtheorem{lemma}{Lemma}%[section]
\newtheorem{corollary}{Corollary}%[section]
\newtheorem{example}{Example}%[section]

\title{\Large Asymptotic Theory for Combining Dependent $p$-Values for Global Hypothesis Testing }
\author{Haoyi Yang and Lingzhou Xue \\ Department of Statistics, The Pennsylvania State University}
\date{}

\begin{document}

\maketitle
 
\begin{abstract}
Combining $p$-values is a fundamental procedure in global hypothesis testing. In modern high-dimensional settings, however, component $p$-values often exhibit complex dependence and rely on asymptotic approximations rather than exact finite-sample uniform distributions. This paper establishes a unified asymptotic theory for weighted transformation statistics that decouples marginal finite-sample approximation error from the joint dependence structure. We also provide sufficient conditions based on conditional probability bounds to verify the joint-tail conditions. Utilizing this framework, we derive explicit dimension-growth and correlation rates for test statistics operating under asymptotic Gaussian and chi-square calibrations. For non-exact finite-sample statistics, we analyze standardized weighted sums, demonstrating how Cramér moderate deviations control relative tail error. Analytical examples demonstrate why both marginal and joint conditions are mathematically indispensable for valid global inference under dependence, and numerical experiments confirm that our asymptotic framework maintains accurate finite-sample size control at extreme significance levels.
\end{abstract}

\noindent\textbf{Keywords:} Combination test; Heavy-tailed transformation; Moderate deviation.

\section{Introduction}

The aggregation of individual $p$-values or test statistics into a global test statistic is a historically important problem in statistical inference. Classical approaches, such as Fisher's, Tippett's, Pearson's, and Stouffer's methods \citep{Fisher1925,tippett1931methods,pearson1934new,stouffer1949adjustment}, alongside the confidence distribution framework \citep{xie2011confidence}, were developed primarily for independent $p$-values. In recent years, the idea of combining asymptotically independent $p$-values or test statistics has drawn considerable attention in high-dimensional hypothesis testing \citep{li2015joint,li2018applications,li2025power,he2021asymptotically,wang2023computationally,yu2020fisher,yu2022power,yu2024power}. A common thread in this literature is the development of power-enhanced tests based on Fisher's method to combine these asymptotically independent $p$-values derived from various high-dimensional tests.

However, individual tests frequently exhibit complex dependence structures in real-world applications, particularly when they are calculated from a shared dataset. For instance, \cite{liu2019acat} pointed out that the variant-level $p$-values from the sequence kernel association test (SKAT) would display moderate or strong correlations. Similarly, \cite{wilson2019harmonic} noted the dependence among the likelihood ratio tests of association between neuroticism and variants, and \cite{liu2020cauchy} showed that the correlation matrix of individual Cochran--Armitage trend tests for the association between the disease status and individual variants would contain very strong correlations. 

To account for dependent tests, classical procedures such as Fisher's method or the minimum $p$-value approach typically require knowledge of the joint null distribution or computationally intensive resampling schemes \citep{yu2025unified}. Furthermore, employing conservative calibrations often leads to a loss of power. To address this, several computationally efficient combination methods have recently emerged in the statistics literature. \cite{liu2020cauchy} proposed the Cauchy combination test (CCT) methods for evaluating the global null hypothesis. 
\cite{wilson2019harmonic} studied the harmonic mean $p$-value (HMP) approach, building on the weighted harmonic mean of individual $p$-values \citep{good1958significance}, as a powerful alternative to Bonferroni correction. \cite{vovk2020combining} and \cite{chen2020trade} showed that valid $p$-values can be combined by scaling up their generalized means to control the familywise error rate (FWER) of multiple dependent tests,
a framework that includes the HMP and CCT as special cases. More recently,
\cite{liu2025heavily} introduced the Half-Cauchy combination test (HCCT) by truncating the left half of the Cauchy distribution. A primary advantage of these statistics is that they can be computed without estimating the full dependency structure among individual $p$-values, although their validity relies on specific assumptions.

Existing literature has established important calibration results for combining exact marginal $p$-values under various dependence assumptions. Assuming pairwise bivariate normality among the individual test statistics, \cite{liu2020cauchy} proved an asymptotic Cauchy tail approximation under the null hypothesis when the number of $p$-values is fixed. Under this same normality assumption and a fixed number of tests, \citet{fang2021heavy} derived tail approximations for combination tests based on regularly varying heavy-tailed distributions. Expanding on this, \cite{gui2025aggregating} studied the combination tests of \cite{fang2021heavy} and established tail
approximations under quasi-asymptotic independence. Concurrently, \cite{long2023cauchy} introduced a theoretical framework for the Cauchy tail
approximation based on the joint lower tails of the exact $p$-values, and \cite{liu2025heavily} subsequently adopted this framework to study their proposed HCCT. On the other hand, \cite{vovk2020combining} and \cite{chen2020trade} showed that suitably calibrated generalized means of valid $p$-values can control the FWER under arbitrary dependence, which is equivalent to testing the global null hypothesis. Unlike the aforementioned works, their theoretical results do not characterize the asymptotic tail behavior of the generalized mean. 

However, these existing theoretical results in the literature do not account for the finite-sample approximation error that arises when component $p$-values are computed using a limiting or approximate reference distribution. Notably, the variant-level statistics from the SKAT in \cite{liu2019acat} and the Cochran–Armitage trend test in \cite{liu2020cauchy} do not satisfy the exact bivariate normality condition; rather, they rely on asymptotic normal approximations using the central limit theorem. In practice, the global test may use significance levels that decrease with the sample size $n$ or the number of tests $d$. At such levels, an approximation error that vanishes in absolute value can remain large relative to the probability being approximated. Thus, marginal weak convergence, or a conventional Berry--Esseen bound stated in absolute error, does not by itself justify replacing a finite-sample test statistic with its limiting reference distribution inside the combination statistic.

Let $h: (0, 1] \to \mathbb{R}$ be a decreasing function such that $h(p)\to\infty$ as $p\downarrow0$. We study weighted combination statistics of the form
\[
T_n=\sum_{i=1}^d w_i h(P_{n,i}),
\]
where $w_i\ge 0$, $\sum_{i=1}^d w_i=1$, and $P_{n,i}$ is the $p$-value of the $i$-th component test. We establish conditions under which
\begin{equation}\label{eq:tail-sum-intro}
 \Pr(T_n>t_n)\approx\sum_{i=1}^d  \Pr\{w_i h(P_{n,i})>t_n\}
\end{equation}
holds, allowing the null tail probability of the combined statistic $T_n$ to be calculated directly from its marginal tail probabilities. Incorporating dependence, however, introduces a separate difficulty. Even with accurate marginal approximations, simultaneous extreme transformed values can occur, which may alter the first-order tail probability of their sum. Moreover, for signed transformations such as the CCT, a large negative value can offset a large positive contribution. Hence, both the marginal approximation error and the relevant joint-tail probabilities must be controlled over the same sequence of diverging thresholds. A crucial insight of our framework is the deliberate decoupling of these requirements: marginal finite-sample calibration on one hand, and the joint dependence structure on the other. By formulating marginal relative error control and pairwise tail control as separate conditions, our approach integrates the dependence analyses typically used for exact $p$-values, while clarifying the additional theory needed to justify the use of approximate component $p$-values.

Our first contribution (Theorem \ref{thm:gentail}) establishes a general tail approximation for the CCT, HCCT, HMP, and generalized means, applicable to both fixed and diverging dimension $d$. Unlike prior frameworks that assume exact finite-sample uniformity, we quantify how finite-sample errors in the component $p$-values and their joint tail probabilities affect the aggregate statistic, providing conditions under which the tail approximation remains valid. The key theoretical insight of this theorem is the decoupling of uniform relative marginal-tail accuracy from pairwise absolute upper-tail control. Furthermore, through two illustrative examples (i.e., Examples 1 and 2 in Section \ref{mainresults}), we demonstrate that these conditions address distinct theoretical issues: marginal misspecification can alter the decay rate of the tail probability under complete independence, whereas extreme dependence can distort the leading constant when all marginal $p$-values are exactly uniform.

Beyond marginal calibration, verifying dependence conditions in Theorem \ref{thm:gentail} poses another bottleneck in the literature. Prior works either assume pairwise bivariate normality to achieve tractability \citep{liu2020cauchy,fang2021heavy} or impose other conditions directly on the joint tails \citep{gui2025aggregating,long2023cauchy,liu2025heavily}. Recognizing that the latter requirement is often analytically intractable, as probabilities of simultaneous extreme events rarely admit simple closed forms, we turn to conditional distributions. In the Gaussian and correlated chi-square reference models studied below, explicit conditional distributions make the joint-tail bounds easier to verify. Motivated by this, our second contribution (Theorem~\ref{thm:gencond}) introduces a novel theoretical framework to verify these dependence structures. Specifically, we establish sufficient conditions based on conditional probability bounds that imply the required joint-tail condition. The core methodological innovation in this theorem is a localization strategy: because large conditioning values can distort the conditional law, we must localize the conditioning event before applying the conditional bound. By bounding the conditional probability uniformly over this localized region, our theory avoids the direct joint-tail evaluation and provides a tractable pathway for proving finite-sample validity under complex dependence.

Our third contribution shows the applicability of our theory across several important classes of statistics, including asymptotically normal or chi-square statistics and weighted sums. To establish these results, the required marginal and joint tail approximations must remain tightly controlled at small probability levels, a regime where standard Berry--Esseen bounds are inadequate due to their generally restrictive rates. For weighted sums, we leverage a Cramér moderate-deviation result to establish the required marginal tail approximation. 

Moreover, our comprehensive numerical experiments demonstrate the finite-sample properties. By examining how this approximation scales with dimension, tail level, and dependence strength, these numerical results complement our theory to illustrate the distinct errors induced by finite-sample marginal calibration and joint-tail behavior.

The rest of this paper is organized as follows. Section~2 presents the notation and preliminaries. Section~3 states our main results, including the general tail-approximation theorem, the sufficient condition based on conditional probabilities, and two examples illustrating why both marginal and pairwise conditions are necessary. Section~4 applies the general theory to specific statistical settings, while Section~5 presents numerical experiments. Section~6 includes a few concluding remarks. The complete proofs are presented in the supplement. 

\section{Preliminaries}

For each $n$, let $S_n=(S_{n,1},\ldots,S_{n,d})^\top$ denote a random vector of component test statistics, where the dimension $d=d_n$ may be fixed or diverging with $n$. Let $F_i$ be the continuous reference cumulative distribution function used to calibrate the $i$-th statistic. The corresponding one-sided reference $p$-value can be written as $P_i(s)=1-F_i(s)$. For notational brevity, we write $P_{n,i}=P_i(S_{n,i})$, where the reference distribution may not coincide with the true finite-sample distribution of $S_{n,i}$. 

Let $h:(0,1)\to\mathbb R$ be a decreasing transformation, and let $w_i>0$, $i=1,\ldots,d$, be positive weights such that $w_i\in(c/d,C/d)$ for fixed constants $0<c<C<\infty$, and $\sum_{i=1}^d w_i=1$. Let $X_{n,i} =h\{P_i(S_{n,i})\}$, and we consider the weighted combination statistic of the form:
\[
T_n:=\sum_{i=1}^{d}w_iX_{n,i} =\sum_{i=1}^{d}w_i h\{P_i(S_{n,i})\}.
\]
We will formally introduce the requirements on the upper-tail inverse of $h$, boundary regularity, and weight allocations in Condition~(A0) of Section 3. Although our primary results are presented for one-sided $p$-values, they readily adapt to two-sided alternatives. For a two-sided test, assuming a reference random variable  $U_i\sim F_i$, the two-sided $p$-value is defined as $P_i^{\pm}(s)=\Pr\{|U_i|\ge |s|\}$, the corresponding combination statistic is
\[
T_n^{\pm}=\sum_{i=1}^dw_i h\{P_i^{\pm}(S_{n,i})\}.
\]
%The modifications needed for $T_n^{\pm}$ are collected after the proofs.
Our framework accommodates several widely used transformations. The Cauchy Combination Test (CCT) employs
$$h_C(p)=\tan\!\left\{\pi\left(\frac12-p\right)\right\},$$
which maps a uniform $p$-value to a standard Cauchy random variable. Notably, this transformation is unbounded both above and below. By contrast, the reciprocal transformation
$$h_H(p)=p^{-1},$$ which underlies the Harmonic Mean $p$-value (HMP), is bounded below. More generally, for any $r>0$ we consider the power transformation $$h_r(p)=p^{-r}.$$ The special case $r=1$ recovers the reciprocal transformation, while the generalized mean construction of \citet{vovk2020combining} utilizes closely related power transformations with calibrations designed to control the FWER under arbitrary dependence. The standard half-Cauchy function provides another example of a transformation that is bounded below:
\[
h_{HC}(p)=\cot\!\left(\frac{\pi p}{2}\right)
\sim\frac{2}{\pi p},\quad \text{as }  p\downarrow0.
\]
Note that non-significant tails will have a minor impact compared to extremely significant tails \citep{good1958significance}, and these four transformations {magnify exceptionally small $p$-values.} The Cauchy, reciprocal, and half-Cauchy transformations have tail index 1, whereas $h_r$ has tail index $1/r$. The tail index dictates the normalization of the leading asymptotic term, while the lower bound of $h$ determines whether negative cancellation effects must be controlled.

\section{Main Results}\label{mainresults}

We will first provide a general tail-approximation theorem in Section 3.1 and then present sufficient conditions via conditional probabilities in Section 3.2. 

\subsection{A General Theory}

Throughout this section, the dimension $d=d_n$ may be fixed or diverging, and we consider a diverging tail threshold $t_n\to\infty$ for the combined statistic, alongside a vanishing auxiliary sequence $\delta_n\to 0$ such that $t_n\delta_n\to\infty.$ This auxiliary sequence dictates the relative perturbation allowed in our tail approximations. 

The following condition on the upper-tail inverse of $h$, boundary regularity, and weight allocations is imposed along the sequence $(n,d_n,t_n,\delta_n)$.

\begin{itemize}

\item[(A0)]
The transformation function $h:(0,1)\rightarrow\mathbb{R}$ is continuous, decreasing, and satisfies that $h(p)\rightarrow\infty$ as $p\rightarrow0.$ There exist $p_h\in(0,1)$ and $u_h\in\mathbb R$ such that the restriction of $h$ to $(0,p_h]$ is strictly decreasing onto $[u_h,\infty)$. For $u\ge u_h$, let $q(u)=h^{-1}(u)$ denote this upper-tail inverse. Define 
$Q_n=\sum_{i=1}^dq(t_n/w_i)$. The probability scale associated with an intermediate summand is required to be comparable to the sum of the individual exceedance probabilities:
\begin{equation}\label{scale-comparable}
   \sup_{a\in[1/C,1/c]}q(at_n)=O(Q_n). 
\end{equation}
Assume further that
\[
\max_{1\le i\le d}
\left|
\frac{h^{-1}\{(1\pm\delta_n)t_n/w_i\}}
     {h^{-1}(t_n/w_i)}-1
\right|\longrightarrow0,
\]
where the convergence is required for both the plus and minus signs.

\end{itemize}

The scale-comparability requirement \eqref{scale-comparable} in Condition (A0) is satisfied for the Cauchy and reciprocal transformations. For the power transformation $h(p)=p^{-r}$, it holds for every $r>0$ when $d$ is fixed. When $d$ diverges under the stated weight bounds, it holds for $r\ge1$ due to $Q_n=t_n^{-1/r}\sum_iw_i^{1/r}$ and the weight bounds. Lemma~\ref{lm:standard-a0} formalizes these calculations, including the stability requirement under the slightly perturbed thresholds in (A0).

\begin{lemma}\label{lm:standard-a0}
Let $t_n\to\infty$, $\delta_n\to0$, and let the weights satisfy the bounds specified in (A0), where the dimension $d=d_n$ may be fixed or diverging as $n\to\infty$. The Cauchy, reciprocal, and half-Cauchy transformations satisfy both the perturbation stability and scale-comparability requirements of Condition (A0). These requirements are also satisfied by the power transformation $h(p)=p^{-r}$ for all $r\ge1$. Specifically, $Q_n\sim\frac{1}{\pi t_n}$ for the Cauchy transformation, $Q_n=\frac{1}{t_n}$ for the reciprocal transformation, $Q_n\sim\frac{2}{\pi t_n}$ for the half-Cauchy transformation, and $Q_n=t_n^{-1/r}\sum_{i=1}^dw_i^{1/r}$for the power transformation with $r\ge1$.
\end{lemma}

Fix a constant $K_0>\max\{C,2/c\}$ so that the auxiliary thresholds used in the theoretical development lie within the interval $I_n=\left(\frac{t_n\delta_n}{K_0},K_0dt_n\right)$. 
The choice of $K_0$ depends on the constants bounding the weights and plays no inferential role. 

The following conditions introduce the marginal and pairwise tail requirements:

\begin{itemize}
\item[(A1)]
The marginal transformed tails exhibit uniformly accurate relative error:
\[
\sup_{1\le i\le d}\;
\sup_{u\in I_n}
\left|
\frac{\Pr(X_{n,i}>u)}{q(u)}
-1
\right|
\longrightarrow0.
\]

\end{itemize}

Condition (A1) is formulated in terms of relative error. Because the denominator $q(u)$ tends to zero, a vanishing absolute error could dominate the probability being approximated. 

\begin{itemize}
\item[(A2)] Simultaneous extreme values must be vanishingly rare. Specifically, the joint probability of two components having large absolute magnitudes must be negligible compared to the target probability of the more extreme individual event:
\[
d^2
P\!\left(
X_{n,i}>u,\,
|X_{n,j}|>v
\right)
=
o\!\left(
q(u)\wedge q(v)
\right)
\]
uniformly over ordered pairs $i\ne j$ and $u,v\in I_n$.

\end{itemize}
Condition (A2) ensures that an extreme value of the combined sum must be driven by exactly one large component. By bounding the joint probability of the absolute values $X_{n,i}$ and $\vert{}X_{n,j}\vert{}$, this condition controls all pairwise interactions that could disrupt the standard sum approximation. Specifically, it guards against two large positive components jointly pushing the sum over the threshold, as well as a large negative component dragging the sum down to mask a legitimate positive exceedance. While negative--negative pairs do not directly threaten the upper tail, bounding the absolute joint tails provides a unified requirement that encompasses all problematic pairwise dependencies.

\begin{theorem}[General tail approximation]\label{thm:gentail}
Assume $t_n\to\infty$, $\delta_n\to 0$, and $t_n\delta_n\to\infty.$ If Conditions (A0), (A1), and (A2) hold, where all three conditions are satisfied uniformly when the dimension $d=d_n$ diverges, then as $n\to\infty$, we have
\[
\frac{\Pr(T_n>t_n)}{\sum_{i=1}^d h^{-1}(t_n/w_i)}\longrightarrow1.
\]
\end{theorem}

Theorem~\ref{thm:gentail} advances two trajectories in the existing literature. The first trajectory focuses on precise tail approximations: \citet{liu2020cauchy} and \citet{fang2021heavy} established foundational results for Cauchy and regularly varying transformations under the assumption of pairwise bivariate normality. While subsequent extensions by \citet{gui2025aggregating}, \citet{long2023cauchy}, and \citet{liu2025heavily} accommodated broader dependence structures by imposing conditions on the joint tails, these frameworks remain restricted to exact marginal $p$-values. In a second trajectory, \citet{vovk2020combining} and \citet{chen2020trade} studied universally valid thresholds rather than characterizing the asymptotic tail behavior of the combination statistic. Theorem~\ref{thm:gentail} characterizes the exact asymptotic tail equivalence while accommodating component $p$-values derived from approximate reference distributions.

In what follows, we provide two examples to illustrate why both marginal and pairwise assumptions (A1) and (A2) address distinct theoretical issues. Example 1 demonstrates that complete independence cannot repair the asymptotic distortions caused by a misspecified marginal tail. On the other hand,  Example 2 maintains uniform marginal $p$-values but does not guarantee the correct leading constant in the presence of extreme tail dependence.

\begin{example}[Misspecified marginal calibration]
Let $S_1,\ldots,S_d$ be independent random variables with tails given by
$\mathbb{P}(S_i > x) = x^{-\alpha}$, for $x \ge 1$ and some $\alpha>0$. Suppose that we incorrectly calibrate these statistics using a standard normal reference distribution, setting $P_i=1-\Phi(S_i)$, and apply the Cauchy transformation $\qquad
h_C(p)=\tan\{\pi(1/2-p)\}$. Since $h^{-1}_C(t)\sim(\pi t)^{-1}$ and
$\Phi^{-1}\{1-h^{-1}_C(t)\}\sim\sqrt{2\log t}$, it follows that
\[
\Pr\{h_C(P_i)>t\}
=\Pr\!\left[S_i>\Phi^{-1}\{1-h^{-1}_C(t)\}\right]
\asymp(2\log t)^{-\alpha/2}.
\]
Thus, even under complete independence between coordinates, the standard denominator in Theorem~\ref{thm:gentail} fails to approximate the tail of $T_n$.
\end{example}

\begin{example}[Complete dependence with a non-unit tail index]
Let $P_1=\cdots=P_d=U$, where $U\sim\operatorname{Unif}(0,1)$, and apply the power transformation $h_r(p)=p^{-r}$ for some $r\ne1$. Because the weights sum to one, the combination statistic simplifies to
\[
T=\sum_{i=1}^dw_iU^{-r}=U^{-r},
\qquad \text{and }\quad
\Pr(T>t)=t^{-1/r}.
\]
By contrast, evaluating the standard tail approximation yields
$\sum_{i=1}^dh^{-1}(t/w_i)
=t^{-1/r}\sum_{i=1}^dw_i^{1/r}$. 
Consequently, the ratio of the true tail to the approximation is
\[
\frac{\Pr(T>t)}{\sum_{i=1}^dh^{-1}(t/w_i)}
=\left(\sum_{i=1}^dw_i^{1/r}\right)^{-1},
\]
which generally differs from 1 unless $r=1$. The pairwise-negligibility condition is therefore essential for transformations with tail index different from one.
\end{example}

When the component statistics exactly follow their continuous reference distributions, the marginal calibration requirement inherently vanishes, yielding the following corollary.

\begin{corollary}[Exact marginal calibration]\label{thm:genbi}
Let $d=d_t$ be fixed or diverge as $t\to\infty$, and let $\delta=\delta_t\to0$ such that $t\delta_t\to\infty$. Suppose that
$S=(S_1,\ldots,S_d)^\top$ has continuous marginal null distribution
functions $F_1,\ldots,F_d$. Define $P_i(S_i)=1-F_i(S_i),$ and $T=\sum_{i=1}^d w_i h\{P_i(S_i)\}.$ 
If Conditions (A0) and (A2) hold with $t_n$ and $\delta_n$
replaced by $t$ and $\delta_t$, respectively, and with both conditions satisfied uniformly if the dimension $d$ diverges, then as $t\to\infty$,
\[
\frac{\Pr(T>t)}
{\sum_{i=1}^d h^{-1}(t/w_i)}
\longrightarrow1.
\]
\end{corollary}

Since each $F_i$ represents the exact continuous null distribution of $S_i$,
the marginal $p$-values $P_i(S_i)$ are uniformly distributed on
$(0,1)$, implying that Condition (A1) is satisfied automatically. Because no finite-sample approximation is involved, the index $n$ is omitted, allowing the result to be formulated directly as an extreme-value limit as $t\to\infty$. Notably, exact Gaussian models are a special case satisfying Conditions~(A0) and (A2). For a fixed dimension $d$, specializing this corollary to exact Gaussian $p$-values under nondegenerate pairwise bivariate
normality recovers the structural form of the tail approximations established by \citet{liu2020cauchy} for the Cauchy transformation and by
\citet{fang2021heavy} for regularly varying transformations. Unlike previous results that rely heavily on the specific properties of bivariate normal distributions, the present corollary extracts the essential dependence requirement into Condition~(A2), thereby accommodating any continuous component statistics.

\subsection{Sufficient Conditions via Conditional Probabilities}

While Theorem~\ref{thm:gentail} establishes the general mechanics of the tail approximation, verifying the joint probability bound in Condition~(A2) is often non-trivial for complex dependence structures. To render these requirements practically verifiable, we now introduce a framework of sufficient conditions formulated via conditional probabilities. Rather than bounding global joint distributions, this approach relies on a localization strategy: by conditioning on a single valid exceedance, we constrain the necessary probability bounds to the relevant auxiliary tail band $I_n$. As established in the forthcoming Theorem~\ref{thm:gencond}, this conditional framework decouples the verification of marginal calibration from the control of pairwise dependence.

Consider a reference random vector $U=(U_1,\ldots,U_d)$ with marginal distribution functions $F_i$ and specified pairwise joint distributions. Define the threshold regions
\[
R_{i,+}(u)=\{s:h(1-F_i(s))>u\},
\qquad \text{and} \quad
R_{i,-}(u)=\{s:h(1-F_i(s))<-u\},
\]
as well as their corresponding marginal tail probabilities
$q_{i,\sigma}(u)=\Pr\{U_i\in R_{i,\sigma}(u)\}$ for $\sigma\in\{+,-\}$. Note that if the transformation $h$ is bounded from below, the set $R_{i,-}(u)$ is empty for all sufficiently large $u$. Under marginal calibration, this reduces to $q_{i,+}(u)=q(u)$.

Condition~(A3) regulates the behavior of the reference vector $U$ and its relationship to the component statistics $S_{n,i}$. Specifically, to achieve this, sub-condition (a) bounds the lower tail, sub-condition (b) enforces an asymptotic conditional independence by partitioning the extreme event space into a typical core and a negligible complement, and sub-condition (c) ensures the finite-sample joint distributions are reliably dominated by the reference model.

\begin{itemize}
\item[(A3)] The following conditions hold uniformly over all ordered pairs $i\ne j$, auxiliary thresholds $u,v\in I_n$, and signs $\sigma\in\{+,-\}$, provided the relevant tail sets are nonempty.
\begin{itemize}
\item[(a)] The lower tail is dominated by the upper tail scale with a fixed constant $C_h$:
\[
q_{i,-}(u)\le C_hq(u).
\]

\item[(b)] Consider the two reference events $\{U_i\in R_{i,+}(u)\}$ and $\{U_j\in R_{j,\sigma}(v)\}$. Let $\{U_k\in R\}$ denote the less likely of the two events and $\{U_\ell\in R'\}$ the other, resolving ties arbitrarily. There exists a regular conditional distribution of $U_\ell$ given $U_k=s$ and a subset $C_k(R)\subseteq R$ such that the probability mass outside $C_k(R)$ is negligible:
\[
\frac{\Pr\{U_k\in R\setminus C_k(R)\}}
     {\Pr(U_k\in R)}=o(d^{-2})
\]
and the conditional probability of the complementary event, given any realization within the core, is uniformly bounded:
\[
\sup_{s\in C_k(R)}
\Pr(U_\ell\in R'\mid U_k=s)=o(d^{-2}).
\]

\item[(c)] The joint tail probabilities are dominated by their corresponding reference probabilities up to a fixed constant $C$:
\[
\Pr\{S_{n,i}\in R_{i,+}(u),S_{n,j}\in R_{j,\sigma}(v)\}
\le C
\Pr\{U_i\in R_{i,+}(u),U_j\in R_{j,\sigma}(v)\}.
\]
The constants $C_h$ and $C$ do not depend on $n,d,i,j,u,$ or $v$.
\end{itemize}
\end{itemize}

For a light-tailed reference distribution, a typical choice for the core subset is $C_k(R)=(b,\kappa b]$ when the threshold region is $R=(b,\infty)$ for some constant $\kappa>1$. The first probability bound in sub-condition (b) requires $\Pr(U_k>\kappa b)$ to be negligible relative to $\Pr(U_k>b)$. The conditional probability in the second bound is then controlled uniformly over the localized interval $b<U_k\le\kappa b$.  As the dimension $d$ diverges, both errors must be $o(d^{-2})$ to ensure that their sum over all $O(d^2)$ ordered pairs remains asymptotically negligible.

The three components of Condition~(A3) can be verified using analytical tools. For instance, sub-condition (a) balances the lower and upper transformed marginal tails,  sub-condition (b) evaluates the conditional reference distribution strictly over a localized core, safely excluding more extreme values that carry negligible probability, and sub-condition (c) requires the joint tail probabilities to be bounded by those of the reference model up to a constant factor. For exact Gaussian and chi-square models, such as those analyzed in Section 4, sub-condition (c) is satisfied naturally. For statistics that rely on asymptotic approximations, verifying sub-condition~(c) requires a separate bivariate tail approximation.

\begin{theorem}[Tail approximation via  conditional probabilities]\label{thm:gencond}
Assume $t_n\to\infty$, $\delta_n\to 0$, and $t_n\delta_n\to\infty$. If Conditions~(A0), (A1), and (A3) hold, with all three conditions satisfied uniformly as the dimension $d=d_n$ diverges, then as $n\to\infty$,
\[
\frac{\Pr(T_n>t_n)}{\sum_{i=1}^d h^{-1}(t_n/w_i)}\longrightarrow1.
\]
\end{theorem}

To illustrate Condition~(A3), we consider two representative reference models. Suppose $(U_i,U_j)$ follows a standard bivariate normal distribution with a uniform correlation bound $\sup_{i\ne j}|\rho_{ij}|\le\rho_0<1.$ The conditional distribution is exactly
\[
U_j\mid U_i=s\sim \mathcal{N}(\rho_{ij}s,1-\rho_{ij}^2).
\]
Fix a constant $\kappa>1$ such that $\rho_0\kappa<1$. Provided the relevant positive and negative normal quantiles are uniformly comparable, restricting the conditioning variable to the localized core $x<U_i\le\kappa x$ leaves a fixed linear gap between the conditional mean and either extreme threshold. Standard Gaussian tail inequalities tightly control this conditional probability, while the marginal tail ratio
\[
\frac{\bar\Phi(\kappa x)}{\bar\Phi(x)}
=\exp\{-\tfrac12(\kappa^2-1)x^2+o(x^2)\}
\]
shows that the probability mass outside the core (i.e., $\Pr(U_i>\kappa x)$) is asymptotically negligible relative to $\Pr(U_i>x)$. Hence, Condition~(A3) is satisfied provided these exponential decay rates dominate $2\log d$.

Similarly, consider a reference model of correlated chi-square statistics: $U_i=\sum_{k=1}^{\nu}X_k^2$, $U_j=\sum_{k=1}^{\nu}Y_k^2,$ and $\operatorname{Corr}(X_k,Y_k)=\rho_{ij}.$ The conditional distribution is a scaled noncentral chi-square:
\[
\frac{U_j}{1-\rho_{ij}^2}\ \bigg|\ U_i=s
\sim
\chi_\nu^2\!\left(\frac{\rho_{ij}^2s}{1-\rho_{ij}^2}\right).
\]
By fixing $\kappa>1$ such that $\rho_0^2\kappa<1$, noncentral chi-square concentration inequalities bound the conditional tail uniformly for the localized core $x<U_i\le\kappa x$. Furthermore, the ratio
\[
\frac{\Pr(U_i>\kappa x)}{\Pr(U_i>x)}
=\exp\{-\tfrac12(\kappa-1)x+o(x)\}
\]
again ensures that the complementary event $\Pr(U_i>\kappa x)$ is negligible relative to $\Pr(U_i>x)$. When the transformation $h$ is unbounded from below, the lower transformed tail corresponds to $U_j$ approaching zero; in this regime, standard noncentral chi-square small-ball bounds verify the tail-balancing requirement of Condition~(A3). Thus, both the Gaussian and correlated chi-square models satisfy the conditional criterion, provided explicit rate conditions link the extreme tail threshold with the diverging dimension $d$.

\section{Explicit Rates under Asymptotic Calibration}

We study two important frameworks in high-dimensional inference: test statistics under asymptotic Gaussian and chi-square calibration settings, and weighted sums operating under asymptotic Gaussian calibration. By deriving explicit rate conditions for these regimes, we demonstrate how our framework accommodates non-exact finite-sample statistics.

\subsection{Asymptotic Gaussian and Chi-Square Calibration}

Gaussian and chi-square approximations are common for deriving component
$p$-values that deviate from exact uniformity in finite samples. We analyze these two fundamental cases, leveraging their conditional reference distributions to verify the dependence requirements of Theorem~\ref{thm:gencond}. Throughout this analysis, Conditions~(A1) and (A3)(c) retain their original roles to govern the marginal approximation error and the comparison of finite-sample joint probabilities and the reference probabilities, respectively.

We assume that the dimension $d\ge2$, the weights $w_i$ satisfy the normalization and bounds established in Section~2, and the transformation function $h$ is one of $h_C$, $h_{HC}$, $h_H$, or $h_r$. For the power transformation $h_r$, $r\ge1$ is
fixed when $d$ diverges, whereas any fixed $r>0$ is allowed for fixed $d$. In the rate conditions that follow, we implicitly set $r=1$ for $h_C$, $h_{HC}$, and $h_H$.
We construct the localized band $I_n$ using the perturbation parameter $\delta_n=1/\log t_n$. To ensure the combination statistic $T_n$ remains finite under the Cauchy transformation $h_C$, we assume that the component $p$-values lie in $(0,1)$ almost surely. 

To govern the interplay between dimension growth and tail dependence, let $0\le\rho_0<1$ be fixed, and suppose
that $t_n\to\infty$. We assume that for some fixed exponent $\eta\ge 0$,
\begin{equation}\label{eq:asymptotic-conditional-rate}
 d=O(t_n^\eta),
 \qquad\text{and }\quad
 \rho_0\sqrt{1+(1+2r)\eta}
 +\sqrt{2r(1-\rho_0^2)\eta}<1.
\end{equation}
The inequality in \eqref{eq:asymptotic-conditional-rate} provides a sufficient range for diverging $d$, while the fixed $d$ regime is recovered by setting $\eta=0$.

For asymptotic Gaussian calibration, let $F_i=\Phi$, and define the reference vector
$U_n$ such that it has standard normal marginals, with every pair jointly normal and correlations satisfying $\max_{i\ne j}|\rho_{n,ij}|\le\rho_0$. The conditional distribution for any pair then becomes
\[
 U_{n,j}\mid U_{n,i}=s
 \sim \mathcal{N}\!\left(\rho_{n,ij}s,1-\rho_{n,ij}^2\right),
\]
which verifies Condition~(A3)(b) under 
the rate condition \eqref{eq:asymptotic-conditional-rate}. The remaining reference
condition (A3)(a) and the transformation condition (A0) are satisfied by the
listed transformations. Thus, the conditional probability requirement
is deduced from the Gaussian reference structure, rather than being
imposed as an additional assumption. For the Cauchy transformation,
positive and negative transformed values correspond to the two tails
of the normal statistic, which are treated by the same conditional
normal calculation. This conditioning structure establishes the tail approximation for asymptotically Gaussian statistics.

\begin{corollary}[Asymptotically Gaussian statistics]
\label{cor:asymptotic-gaussian}
Let $P_{n,i}=1-\Phi(S_{n,i})$. If Conditions (A1) and (A3)(c) hold
with respect to the Gaussian reference vector $U_n$ defined above, and the rate condition
\eqref{eq:asymptotic-conditional-rate} is satisfied, then as $n\to\infty$,
\begin{equation*}\label{eq:asymptotic-combination-tail}
 \frac{\Pr(T_n>t_n)}
 {\sum_{i=1}^d h^{-1}(t_n/w_i)}\longrightarrow1.
\end{equation*}
\end{corollary}

For asymptotic chi-square calibration, fix an integer $\nu\ge1$. Let
$Y_{n,1},\ldots,Y_{n,d}$ be $\nu$-dimensional standard Gaussian
vectors with the following pairwise covariance structure:
\[
 \begin{pmatrix}Y_{n,i}\\Y_{n,j}\end{pmatrix}
 \sim \mathcal{N}\!\left(0,
 \begin{pmatrix}
 I_\nu&\rho_{n,ij}I_\nu\\
 \rho_{n,ij}I_\nu&I_\nu
 \end{pmatrix}\right),
 \qquad \max_{i\ne j}|\rho_{n,ij}|\le\rho_0.
\]
Define the squared-norm reference statistics $U_{n,i}=\|Y_{n,i}\|^2$ with marginal distributions $F_i=F_{\chi_\nu^2}$. In this case, rotational invariance results in the noncentral chi-square conditional distribution with the noncentrality parameter $\frac{\rho_{n,ij}^2s}{1-\rho_{n,ij}^2}$:
\[
 \frac{U_{n,j}}{1-\rho_{n,ij}^2}
 \;\bigg|\; U_{n,i}=s
 \sim\chi_\nu^2\!\left(
 \frac{\rho_{n,ij}^2s}{1-\rho_{n,ij}^2}\right),
\]
which verifies Condition~(A3)(b). Unlike the Gaussian calibration setting, a large negative Cauchy-transformed
value now corresponds to a statistic localized near zero. This does not
require an additional assumption on the chi-square reference law.
To see this, write $F=F_{\chi_\nu^2}$. For each fixed $a_0>0$, there
exists a constant $K$, depending only on $\nu$, $\rho_0$, and $a_0$, such
that, uniformly in $n$ and for $i\ne j$, the lower-tail and cross-tail joint probabilities satisfy
\begin{align*}
 \Pr(U_{n,i}\le a,U_{n,j}\le b)
 &\le K F(a)F(b),\\
 \Pr(U_{n,i}\le a,U_{n,j}>x)
 &\le K F(a)\{1-F(x)\},
\end{align*}
for $0<a,b\le a_0$ and $x\ge0$. These bounds follow from the
bounded conditional density of the underlying Gaussian vector and the  property that 
$F(a)\asymp a^{\nu/2}$ near zero. They remain valid for $\nu=1$, despite the marginal chi-square density itself diverging to infinity at zero. For
non-exact statistics, Condition (A3)(c) transfers these joint
probability bounds from the continuous reference model to the finite-sample distributions. This noncentral conditioning structure establishes the tail approximation for asymptotically chi-square statistics.

\begin{corollary}[Asymptotically chi-square statistics]
\label{cor:asymptotic-chisquare}
Let $P_{n,i}=1-F_{\chi_\nu^2}(S_{n,i})$. If Conditions (A1) and
(A3)(c) hold with respect to the chi-square reference vector $U_n$ defined above and
the rate condition \eqref{eq:asymptotic-conditional-rate} is satisfied, then as $n\to\infty$,
\begin{equation*}
 \frac{\Pr(T_n>t_n)}
 {\sum_{i=1}^d h^{-1}(t_n/w_i)}\longrightarrow1.
\end{equation*}
\end{corollary}

\subsection{Weighted Sums with Gaussian Calibration}

We analyze component $p$-values derived via Gaussian approximations rather than exact finite-sample reference distributions. For every $n$, let $\xi_{n,1},\ldots,\xi_{n,n}$ be \emph{i.i.d.} random vectors in $\mathbb R^d$, where $\xi_{n,ij}$ is the $i$-th coordinate of $\xi_{n,j}$. For deterministic weights $b_{n,ij}$, define the weighted sums and their maximum absolute weight as
\[
S_{n,i}=\sum_{j=1}^n b_{n,ij}\xi_{n,ij},
\qquad
\text{and} \quad
L_n=\max_{1\le i\le d}\max_{1\le j\le n}|b_{n,ij}|,
\]
where $\sum_{j=1}^n b_{n,ij}^2=1$. We impose the following structural condition on the array.

\begin{itemize}
\item[(A4)] For every $n$, $1\le i\le d$, and $1\le j\le n$, $E(\xi_{n,ij})=0$ and $\operatorname{Var}(\xi_{n,ij})=1$. There exist fixed constants $\lambda_0,C_0>0$ such that the moment generating functions are uniformly bounded:
\(
\sup_n\max_{1\le i\le d}\max_{1\le j\le n}
E\exp\{\lambda_0|\xi_{n,ij}|\}\le C_0.
\)
Assume that $L_n\to0$.
\end{itemize}

Condition (A4) governs the accuracy of the marginal Gaussian approximation. Dependence among the coordinates of each vector $\xi_{n,j}$ is allowed, but it must be controlled separately via the pairwise joint-tail requirement in Condition~(A2).

\begin{lemma}[Relative moderate-deviation approximation]\label{lm6}
Under Condition~(A4), for every deterministic sequence $r_n\ge0$ satisfying $r_n=o(L_n^{-1/3})$, we have as $n\to\infty$,
\[
\sup_{1\le i\le d}\sup_{0\le u\le r_n}
\left|\frac{\Pr(S_{n,i}>u)}{1-\Phi(u)}-1\right|\longrightarrow0.
\]
\end{lemma}

Lemma~\ref{lm6} provides a weighted triangular-array formulation of the classical Cram\'er moderate-deviation expansion. Its proof is detailed below to show the uniformity over the diverging coordinates and to clarify the critical role of the maximum coefficient $L_n$.

For any $s\in I_n$, define the underlying Gaussian quantile corresponding to the transformed threshold $s$ as
\(
z_n(s)=\Phi^{-1}\{1-q(s)\}.
\)
Define the transformed combination components and the aggregate statistic as 
\[
X_{n,i}=h\{1-\Phi(S_{n,i})\},
\qquad \text{and} \quad
T_n=\sum_{i=1}^d w_iX_{n,i}.
\]
\begin{theorem}[Weighted sums]\label{thm4}
Under the asymptotic regime stated in Section~3, suppose that Conditions (A0) and (A4) hold, and also that 
\(
\sup_{s\in I_n} z_n(s)=o(L_n^{-1/3}).
\)
If these transformed variables satisfy Condition (A2), then as $n\to\infty$,
\[
\frac{\Pr(T_n>t_n)}{\sum_{i=1}^d h^{-1}(t_n/w_i)}\longrightarrow1.
\]
This tail equivalence applies to both fixed and diverging $d$, provided all conditions hold uniformly over the coordinates and ordered pairs.
\end{theorem}

Lemma~\ref{lm6} provides the relative marginal-tail approximation required by Condition~(A1). Classical Berry--Esseen bounds control absolute error and are consequently too coarse when the target tail probability approaches zero. In contrast, Lemma~\ref{lm6} guarantees vanishing relative error over the domain $u=o(L_n^{-1/3})$. When the weights are uniform, $L_n=O(n^{-1/2})$, this range expands to $u=o(n^{1/6})$. Theorem~\ref{thm4} still explicitly assumes Condition~(A2), as marginal moderate-deviation guarantees do not inherently govern joint-tail probabilities.

\section{Numerical Experiments}
We evaluate the empirical performance of five $p$-value combination methods: the Cauchy combination test (CCT), the harmonic mean $p$-value (HMP), Fisher's combination test (FCT), and two generalized mean procedures. Given component $p$-values $(P_1, \ldots, P_d)^\top$, the combined statistics for the HMP and CCT are defined as 
$p_{\rm HMP}=\left(\sum_{i=1}^dw_iP_i^{-1}\right)^{-1}$ and $p_{\rm CCT}=\frac12-\frac1\pi\arctan\!\left[\sum_{i=1}^dw_i
\tan\!\left\{\pi\left(\frac12-P_i\right)\right\}\right]$.
%Here $p_{\rm HMP}$ is the unadjusted harmonic-mean quantity; it should be distinguished from calibrations designed to guarantee validity under arbitrary dependence.
FCT is evaluated using its classical chi-square calibration, which intrinsically assumes marginal independence. For the generalized means, we define
$M_{-r}=\left(\sum_{i=1}^dw_iP_i^{-r}\right)^{-1/r}$  and employ the exponent $r=5$. We consider two calibrations of this statistic: KAP-W, which denotes the dependence-robust calibration of \citet{vovk2020combining}, and KAP-P, which implements $p_{\rm KAP-P}=\min\left\{1,
\left(\sum_{i=1}^dw_i^{1/r}\right)M_{-r}\right\}$, the asymptotic calibration derived from our Theorem~\ref{thm:gencond}. The following tables report empirical rejection proportions at the nominal levels $\alpha \in \{0.05, 0.01, 0.001\}$.

\begin{table}[H]
\centering
\caption{Empirical sizes under the Gaussian model with AR(1) dependence ($\rho=0.9$).}
\label{tab1}
\begin{tabular}{ll ccccc}
\toprule
& & \multicolumn{5}{c}{Dimension $d$} \\
\cmidrule(lr){3-7}
Nominal Level & Method & 100 & 200 & 400 & 600 & 800 \\
\midrule
\multirow{5}{*}{$\alpha=0.05$} 
& CCT   & 0.064 & 0.064 & 0.063 & 0.062 & 0.061 \\
& HMP   & 0.078 & 0.084 & 0.089 & 0.092 & 0.094 \\
& FCT   & 0.293 & 0.307 & 0.319 & 0.325 & 0.326 \\
& KAP-W & 0.021 & 0.022 & 0.022 & 0.023 & 0.024 \\
& KAP-P & 0.025 & 0.027 & 0.027 & 0.029 & 0.029 \\
\midrule
\multirow{5}{*}{$\alpha=0.01$} 
& CCT   & 0.0118 & 0.0119 & 0.0118 & 0.0118 & 0.0116 \\
& HMP   & 0.0124 & 0.0126 & 0.0128 & 0.0129 & 0.0128 \\
& FCT   & 0.237  & 0.253  & 0.264  & 0.270  & 0.269  \\
& KAP-W & 0.0045 & 0.0050 & 0.0050 & 0.0053 & 0.0054 \\
& KAP-P & 0.0055 & 0.0061 & 0.0062 & 0.0064 & 0.0066 \\
\midrule
\multirow{5}{*}{$\alpha=0.001$} 
& CCT   & 0.0011 & 0.0011 & 0.0011 & 0.0011 & 0.0011 \\
& HMP   & 0.0011 & 0.0011 & 0.0011 & 0.0011 & 0.0011 \\
& FCT   & 0.183  & 0.196  & 0.207  & 0.212  & 0.214  \\
& KAP-W & 0.0005 & 0.0005 & 0.0006 & 0.0006 & 0.0005 \\
& KAP-P & 0.0006 & 0.0007 & 0.0008 & 0.0007 & 0.0006 \\
\bottomrule
\end{tabular}
\end{table}

\paragraph{Exact Gaussian calibration.} We generate $d$-dimensional random vectors $X\sim \mathcal{N}(0,\Sigma)$ with the autoregressive covariance structure $\Sigma_{ij}=\rho^{|i-j|}$, and construct one-sided marginal $p$-values $P_i=1-\Phi(X_i)$. We evaluate dimensions $d\in\{100,200,400,600,800\}$ under both strong positive and alternating correlation scenarios, setting $\rho\in\{0.9,-0.9\}$. The empirical results in Tables~\ref{tab1}--\ref{tab2} show that calibration accuracy depends on both the aggregation method and the targeted tail level. Under strong positive dependence ($\rho=0.9$), the CCT converges toward nominal calibration as the threshold $\alpha$ becomes more extreme. The HMP, while noticeably liberal at $\alpha=0.05$, mirrors the CCT behavior at $\alpha=0.001$. Conversely, under alternating correlation ($\rho=-0.9$), the induced negative dependence structure causes the CCT to become highly conservative. As anticipated, Fisher's FCT, which relies strictly on an independence calibration, exhibits severe anti-conservative bias across both dependence structures. The KAP-W and KAP-P procedures remain universally conservative; this behavior is structurally expected, as they are engineered to guarantee validity across broad, arbitrary dependence classes rather than achieving exact asymptotic tail equivalence.

\begin{table}[H]
\centering
\caption{Empirical sizes under the Gaussian model with AR(1) dependence ($\rho=-0.9$).}
\label{tab2}
\begin{tabular}{ll ccccc}
\toprule
& & \multicolumn{5}{c}{Dimension $d$} \\
\cmidrule(lr){3-7}
Nominal Level & Method & 100 & 200 & 400 & 600 & 800 \\
\midrule
\multirow{5}{*}{$\alpha=0.05$} 
& CCT   & 0.017 & 0.018 & 0.019 & 0.020 & 0.021 \\
& HMP   & 0.078 & 0.083 & 0.088 & 0.092 & 0.093 \\
& FCT   & 0.100 & 0.104 & 0.106 & 0.108 & 0.108 \\
& KAP-W & 0.029 & 0.030 & 0.030 & 0.031 & 0.031 \\
& KAP-P & 0.035 & 0.036 & 0.037 & 0.038 & 0.038 \\
\midrule
\multirow{5}{*}{$\alpha=0.01$} 
& CCT   & 0.0041 & 0.0041 & 0.0042 & 0.0046 & 0.0048 \\
& HMP   & 0.0124 & 0.0116 & 0.0117 & 0.0122 & 0.0121 \\
& FCT   & 0.0479 & 0.0469 & 0.0462 & 0.0464 & 0.0453 \\
& KAP-W & 0.0063 & 0.0061 & 0.0063 & 0.0067 & 0.0068 \\
& KAP-P & 0.0080 & 0.0075 & 0.0077 & 0.0082 & 0.0084 \\
\midrule
\multirow{5}{*}{$\alpha=0.001$} 
& CCT   & 0.0004 & 0.0005 & 0.0005 & 0.0005 & 0.0006 \\
& HMP   & 0.0010 & 0.0011 & 0.0011 & 0.0010 & 0.0011 \\
& FCT   & 0.0195 & 0.0185 & 0.0166 & 0.0157 & 0.0148 \\
& KAP-W & 0.0007 & 0.0007 & 0.0007 & 0.0007 & 0.0007 \\
& KAP-P & 0.0008 & 0.0009 & 0.0008 & 0.0008 & 0.0010 \\
\bottomrule
\end{tabular}
\end{table}

\paragraph{Weighted Bernoulli sums.} Transitioning from exact reference models, we now evaluate the performance of combination tests under the asymptotic Gaussian calibration regime established in Section~4.2. We generate $n=400$ independent observations of correlated Bernoulli vectors $X_j \in \{0,1\}^d$, fixing the marginal success probability at $p=1/2$ and the pairwise within-subject correlation at $\rho=0.5$. For each coordinate $i$, we compute the standardized sum and its corresponding asymptotic marginal $p$-value: $S_i=\frac{\sum_{j=1}^nX_{ij}-n/2}{\sqrt{n/4}}$ and $P_i=1-\Phi(S_i)$. We evaluate the same dimensional sequence $d\in\{100, 200, 400, 600, 800\}$. Because these component $p$-values rely on a central limit approximation, they deviate from exact uniformity in finite samples. Consequently, the empirical rejection proportions reported in Table~\ref{tab3} reflect the compounded effects of marginal finite-sample approximation error and high-dimensional joint dependence. Consistent with Theorem~\ref{thm4}, while both the CCT and HMP exhibit noticeable size inflation at $\alpha=0.05$, their empirical sizes align closely with the nominal level at $\alpha=0.001$, confirming the accuracy of the asymptotic approximations.

\begin{table}[H]
\centering
\caption{Empirical rejection proportions for weighted Bernoulli sums with pairwise $\rho=0.5$.}
\label{tab3}
\begin{tabular}{ll ccccc}
\toprule
& & \multicolumn{5}{c}{Dimension $d$} \\
\cmidrule(lr){3-7}
Nominal Level & Method & 100 & 200 & 400 & 600 & 800 \\
\midrule
\multirow{5}{*}{$\alpha=0.05$} 
& CCT   & 0.087 & 0.093 & 0.099 & 0.102 & 0.105 \\
& HMP   & 0.095 & 0.102 & 0.109 & 0.112 & 0.116 \\
& FCT   & 0.307 & 0.341 & 0.366 & 0.378 & 0.386 \\
& KAP-W & 0.033 & 0.031 & 0.027 & 0.026 & 0.033 \\
& KAP-P & 0.045 & 0.042 & 0.038 & 0.036 & 0.033 \\
\midrule
\multirow{5}{*}{$\alpha=0.01$} 
& CCT   & 0.0138 & 0.0149 & 0.0155 & 0.0161 & 0.0169 \\
& HMP   & 0.0141 & 0.0150 & 0.0157 & 0.0163 & 0.0172 \\
& FCT   & 0.256  & 0.303  & 0.340  & 0.356  & 0.366  \\
& KAP-W & 0.0077 & 0.0069 & 0.0057 & 0.0059 & 0.0068 \\
& KAP-P & 0.0078 & 0.0101 & 0.0085 & 0.0080 & 0.0069 \\
\midrule
\multirow{5}{*}{$\alpha=0.001$} 
& CCT   & 0.0012 & 0.0011 & 0.0012 & 0.0012 & 0.0012 \\
& HMP   & 0.0012 & 0.0011 & 0.0012 & 0.0012 & 0.0012 \\
& FCT   & 0.209  & 0.266  & 0.309  & 0.331  & 0.345  \\
& KAP-W & 0.0007 & 0.0008 & 0.0006 & 0.0009 & 0.0008 \\
& KAP-P & 0.0011 & 0.0008 & 0.0010 & 0.0009 & 0.0008 \\
\bottomrule
\end{tabular}
\end{table}

\begin{table}[H]
\centering
\caption{Empirical rejection proportions under the exact Gaussian model with AR(1) dependence for varying correlation parameters $\rho$.}
\label{tab4}
\begin{tabular}{ll ccccc}
\toprule
& & \multicolumn{5}{c}{Autoregressive Parameter $\rho$} \\
\cmidrule(lr){3-7}
Nominal Level & Method & 0.0 & 0.2 & 0.4 & 0.6 & 0.8 \\
\midrule
\multirow{5}{*}{$\alpha=0.05$} 
& CCT   & 0.049 & 0.050 & 0.052 & 0.055 & 0.059 \\
& HMP   & 0.084 & 0.086 & 0.087 & 0.090 & 0.091 \\
& FCT   & 0.049 & 0.083 & 0.129 & 0.188 & 0.268 \\
& KAP-W & 0.039 & 0.039 & 0.039 & 0.038 & 0.031 \\
& KAP-P & 0.048 & 0.049 & 0.049 & 0.047 & 0.038 \\
\midrule
\multirow{5}{*}{$\alpha=0.01$} 
& CCT   & 0.0099 & 0.0104 & 0.0100 & 0.0104 & 0.0112 \\
& HMP   & 0.0111 & 0.0116 & 0.0112 & 0.0117 & 0.0124 \\
& FCT   & 0.0100 & 0.0249 & 0.0548 & 0.1062 & 0.1940 \\
& KAP-W & 0.0079 & 0.0083 & 0.0079 & 0.0078 & 0.0066 \\
& KAP-P & 0.0099 & 0.0103 & 0.0098 & 0.0097 & 0.0082 \\
\midrule
\multirow{5}{*}{$\alpha=0.001$} 
& CCT   & 0.0011 & 0.0010 & 0.0010 & 0.0011 & 0.0011 \\
& HMP   & 0.0011 & 0.0010 & 0.0010 & 0.0011 & 0.0011 \\
& FCT   & 0.0008 & 0.0046 & 0.0173 & 0.0507 & 0.1302 \\
& KAP-W & 0.0008 & 0.0008 & 0.0008 & 0.0008 & 0.0007 \\
& KAP-P & 0.0011 & 0.0010 & 0.0010 & 0.0010 & 0.0009 \\
\bottomrule
\end{tabular}
\end{table}

\paragraph{Effect of dependence strength.}
To investigate the influence of dependence strength, we fix the dimension at $d=400$ and the sample size at $n=400$ for the asymptotic setting while varying the underlying pairwise dependence parameter $\rho\in\{0, 0.2, 0.4, 0.6, 0.8\}$. Tables~\ref{tab4} and~\ref{tab5} document the resulting empirical sizes under the exact Gaussian and standardized Bernoulli frameworks, respectively. As positive dependence increases, the empirical size of FCT inflates severely, demonstrating the inadequacy of independence-based calibrations. In contrast, the CCT remains highly robust. In the Gaussian experiment, its empirical sizes track the nominal levels across the correlation spectrum. Under the Bernoulli data-generating mechanism, the rejection proportions do not scale strictly monotonically with $\rho$. 

\begin{table}[H]
\centering
\caption{Empirical rejection proportions for weighted Bernoulli sums for varying pairwise correlation parameters $\rho$.}
\label{tab5}
\begin{tabular}{ll ccccc}
\toprule
& & \multicolumn{5}{c}{Correlation Parameter $\rho$} \\
\cmidrule(lr){3-7}
Nominal Level & Method & 0.0 & 0.2 & 0.4 & 0.6 & 0.8 \\
\midrule
\multirow{5}{*}{$\alpha=0.05$} 
& CCT   & 0.049 & 0.068 & 0.097 & 0.093 & 0.070 \\
& HMP   & 0.084 & 0.099 & 0.114 & 0.099 & 0.072 \\
& FCT   & 0.050 & 0.308 & 0.362 & 0.366 & 0.355 \\
& KAP-W & 0.034 & 0.033 & 0.031 & 0.022 & 0.010 \\
& KAP-P & 0.050 & 0.049 & 0.043 & 0.030 & 0.013 \\
\midrule
\multirow{5}{*}{$\alpha=0.01$} 
& CCT   & 0.0097 & 0.0107 & 0.0139 & 0.0164 & 0.0141 \\
& HMP   & 0.0100 & 0.0110 & 0.0142 & 0.0165 & 0.0142 \\
& FCT   & 0.0096 & 0.247  & 0.329  & 0.343  & 0.336  \\
& KAP-W & 0.0062 & 0.0065 & 0.0058 & 0.0046 & 0.0025 \\
& KAP-P & 0.0095 & 0.0099 & 0.0091 & 0.0069 & 0.0032 \\
\midrule
\multirow{5}{*}{$\alpha=0.001$} 
& CCT   & 0.0011 & 0.0010 & 0.0009 & 0.0013 & 0.0013 \\
& HMP   & 0.0011 & 0.0010 & 0.0009 & 0.0013 & 0.0013 \\
& FCT   & 0.0010 & 0.189  & 0.294  & 0.321  & 0.322  \\
& KAP-W & 0.0007 & 0.0005 & 0.0005 & 0.0007 & 0.0003 \\
& KAP-P & 0.0011 & 0.0009 & 0.0008 & 0.0009 & 0.0004 \\
\bottomrule
\end{tabular}
\end{table}

\section{Conclusion}

This paper establishes a unified theoretical framework linking component-level calibration and pairwise tail control to the asymptotic null tail of weighted transformation statistics. Our primary theorem strictly decouples finite-sample marginal approximation error from the high-dimensional dependence structure. The conditional-probability framework decomposes the joint-tail requirement into explicitly verifiable geometric steps, while our constructed counterexamples demonstrate why both marginal and joint conditions remain mathematically indispensable. Exact Gaussian and correlated chi-square models instantiate this abstract dependence analysis, yielding explicit constraints on dimension growth and correlation strength. For non-exact statistics, standardized weighted sums illustrate how Cram\'er moderate deviations control the relative error of approximate $p$-values at tail thresholds.

\bibliographystyle{agsm} 
\bibliography{cauchy}

\end{document}